\documentclass[12pt]{article}
\usepackage{amssymb}
\newcommand{\ft}[2]{{\textstyle\frac{#1}{#2}}}
\def\diag{\mathop{\rm diag}\nolimits}

\def\Aut{\mathop{\rm Aut}\nolimits}

\def\rmi{{\rm i}}
\def\rmd{{\rm d}}
\newsavebox{\uuunit}
\sbox{\uuunit}
    {\setlength{\unitlength}{0.825em}
     \begin{picture}(0.6,0.7)
        \thinlines
        \put(0,0){\line(1,0){0.5}}
        \put(0.15,0){\line(0,1){0.7}}
        \put(0.35,0){\line(0,1){0.8}}
       \multiput(0.3,0.8)(-0.04,-0.02){12}{\rule{0.5pt}{0.5pt}}
     \end {picture}}
\newcommand {\unity}{\mathord{\!\usebox{\uuunit}}}
\csname @addtoreset\endcsname{equation}{section}
\newcommand{\SU}{\mathop{\rm SU}}
\newcommand{\SO}{\mathop{\rm SO}}
\newcommand{\U}{\mathop{\rm {}U}}
\newcommand{\USp}{\mathop{\rm {}USp}}

\newcommand{\Sl}{\mathop{\rm {}S}\ell }
\newcommand{\Gl}{\mathop{\rm {}G}\ell }
\newcommand{\Purple}[1]{#1}
\newcommand{\Red}[1]{#1}
\newcommand{\OliveGreen}[1]{#1}
\newcommand{\Magenta}[1]{#1}
\newcommand{\Blue}[1]{#1}
\newcommand{\Maroon}[1]{#1}
\newif\ifpdf
\ifx\pdfoutput\undefined
   \pdffalse
   \usepackage{cite}
  \usepackage{epsf}
 \else
   \pdfoutput=1
   \pdftrue
  \usepackage[pdftex]{hyperref,graphicx}
\fi

\begin{document}

 %%%%%%%%%%%%%%%%%%%%%%%%%%%%%%%%%%%%%%%%%%%%%%%%%%%%%%%%%%%
\begin{titlepage}
\begin{flushright}
UG-05-08\\
KUL-TF-05/27\\
ITP-UU-05/55\\
SPIN-05/35\\
math-dg/0512084
\end{flushright}
\vspace{.5cm}
\begin{center}
\baselineskip=16pt {\LARGE    The identification of conformal
hypercomplex \\ \vskip 0.2cm and quaternionic manifolds
}\\
\vfill%\vskip 15mm%27.mm
{\large Eric Bergshoeff$^1$,  Stefan Vandoren$^2$ %\\[2mm]
and Antoine Van Proeyen$^3$
  } \\
\vskip 5mm
{\small  $^1$ Center for Theoretical Physics, University of Groningen,\\
       Nijenborgh 4, 9747 AG Groningen, The Netherlands. \\ [2mm]
   %    e.bergshoeff@phys.rug.nl,
   $^2$ Institute for Theoretical Physics, Utrecht University, \\
 Leuvenlaan 4, 3508 TA Utrecht, The Netherlands. \\[2mm]
  $^3$ Instituut voor Theoretische Fysica, Katholieke Universiteit Leuven,\\
       Celestijnenlaan 200D B-3001 Leuven, Belgium.
 }
\end{center}
\vfill
\begin{center}
{\bf Abstract}
\end{center}
{\small
 We review the map between  hypercomplex manifolds that admit a
closed homothetic Killing vector (i.e.\ `conformal hypercomplex'
manifolds) and quaternionic manifolds of 1 dimension less. This map is
related to a method for constructing supergravity theories using
superconformal techniques. An explicit relation between the structure of
these manifolds is presented, including curvatures and symmetries. An
important role is played by `$\xi$ transformations', relating connections
on quaternionic manifolds, and a new type `$\hat\xi$ transformations'
relating complex structures on conformal hypercomplex manifolds. In this
map, the subclass of conformal hyper-K{\"a}hler manifolds is mapped to
quaternionic-K{\"a}hler manifolds.
 } \vspace{2mm} \vfill \hrule width 3.cm
 \vspace{1mm}

\noindent Contribution to the  proceedings volume for the Conference
"Symmetry in Geometry and Physics" in honour of Dmitri Alekseevsky,
September 2005.

 {\footnotesize \noindent
e-mails:
 e.a.bergshoeff@rug.nl,  s.vandoren@phys.uu.nl,
antoine.vanproeyen@fys.kuleuven.be }

\end{titlepage}
\addtocounter{page}{1}
 \tableofcontents{}
\newpage
%%%%%%%%%%%%%%%%
\section{Introduction}
The paper \cite{Alekseevsky1975} of Dmitri Alekseevsky received a lot of
attention in the physics literature as the homogeneous
quaternionic-K{\"a}hler spaces that he investigated occur in $N=2$
supergravity theories \cite{Bagger:1983tt}. Further applications of
quaternionic geometry in supersymmetric theories have been discussed in
\cite{Alekseevsky:2001if}.

In 1996, a general framework to study hypercomplex, hyper-K{\"a}hler,
quaternionic and quaternionic-K{\"a}hler manifolds was given in
\cite{AM1996}. These are commonly denoted by `qua\-ter\-ni\-onic-like
manifolds'. Our paper \cite{Bergshoeff:2004nf} is in many respects a
continuation of \cite{AM1996}. We will review here the main results of
this work.

In this paper, we explain the 1-to-1 correspondence\footnote{As will be
explained below, the correspondence is actually 1-to-1 between families
(or `equivalence classes') of manifolds.} (locally) between conformal
hypercomplex manifolds of quaternionic dimension $n_H+1$ and quaternionic
manifolds of dimension $n_H$. Furthermore, we show that this 1-to-1
correspondence is also applicable between the subset of hypercomplex
manifolds that are hyper-K{\"a}hler and the subset of quaternionic manifolds
that are quaternionic-K{\"a}hler. The map between quaternionic-K{\"a}hler and
hyper-K{\"a}hler manifolds is constructed by Swann \cite{Swann}, and its
generalization to quaternionic manifolds is treated in
\cite{PedersenPS1998}. Here, we give explicit expressions for the complex
structures and connections, curvatures and symmetries. These results are
useful in the context of the conformal tensor calculus in supergravity
\cite{deWit:1985px,Galicki:1992tm,deWit:1999fp,deWit:2001dj,Bergshoeff:2004kh}.

Our work was initiated by an investigation of couplings of
hypermultiplets in supergravity. In section \ref{ss:hypermultiplets} we
will give an overview of the geometric structures in supergravity
theories and indicate where quaternionic geometry finds it place. Section
\ref{ss:quatlike} is for a large part a review of \cite{AM1996}. We give
more details on the vielbeins in these manifolds, which is necessary to
discuss supersymmetry. We will show the need of torsionless affine
connections when discussing supersymmetric theories in spacetime
dimensions $D=3,4,5$ and 6. We will also give results on the relations
between curvature decompositions of the quaternionic-like manifolds.

The main part of this review is section \ref{ss:map}, where the precise
correspondence between conformal hypercomplex manifolds of dimension
$n_H+1$ and quaternionic manifolds of dimension $n_H$ is explained. We
start that section by explaining the relevance of closed homothetic
Killing vectors. Then the general structure of the map is exhibited. We
explain the relevance of the $\xi $ transformations in the quaternionic
manifolds and the existence of similar so-called $\hat\xi$
transformations in conformal hypercomplex manifolds. We finish that
section with a pictorial representation of the map.  A short treatment of
the symmetries of these manifolds is given in section
\ref{ss:symmetries}. Such symmetries are a generalization of isometries
that occur in manifolds with a metric.

We finish in section \ref{ss:conclusions} with conclusions and some
remarks on the relevance of the signature of the extra quaternion in the
hypercomplex manifold.

\section[Hypermultiplets and hypercomplex/quaternionic  manifolds]
{Hypermultiplets and hypercomplex/quaternionic \\ manifolds}

 \label{ss:hypermultiplets}

Table \ref{tbl:mapsusy} gives an overview of theories with rigid and
local supersymmetry\footnote{This table is a shorter version of a table
in \cite{VanProeyen:2003zj}. } in dimensions $D=4$ to $D=11$. The latter
is the maximal dimension for supersymmetric field theories. The top row
indicates the number of real independent components of the spinors
describing the supersymmetry generators. The lowest row indicates which
theories exist only for supergravity, or for supersymmetry and
supergravity. Supergravity is the theory of \textit{local} supersymmetry,
i.e.\ where there is supersymmetry invariance for transformations that
can differ in each spacetime point, as opposed to \textit{rigid}
supersymmetry where the same transformation should be applied for any
point of spacetime. We will concentrate on the theories with 8
supercharges for a reason that we will now explain.

%%%%%%%%%%%%%%%%%%%%%%%%%%%%%%%%%%%%%%%%%%%%%%%%%
\begin{table}[htbp]
  \caption{{\it Supersymmetry and supergravity theories in dimensions 4 to
  11.}  An entry represents the possibility to have supergravity
theories in a specific spacetime dimension $D$ with the number of
supersymmetries indicated in the top row. At the bottom is indicated
whether these theories exist only in supergravity, or also with just
rigid supersymmetry.}
  \label{tbl:mapsusy}
\begin{center}\tabcolsep 5pt
  \begin{tabular}{| *{10}{c|} }
\hline
 $D$  & \multicolumn{2}{c|}{32} & 24  & 20 & \multicolumn{2}{c|}{16}  & 12 & 8 & 4  \\
\hline
11   & M & \multicolumn{1}{c|}{ } & &  & \multicolumn{2}{c|}{ }  &  &  &  \\
10  & IIA & IIB & &  &I&  &  &  &  \\
9  &  \multicolumn{2}{c|}{$N=2$ } &  &  &
  $N=1$&  &  &  &  \\
8   &  \multicolumn{2}{c|}{$N=2$ }&  &  &
$N=1$ &  &  &  &  \\
7   &  \multicolumn{2}{c|}{$N=4$ } &  &  &
$N=2$ & &  &  &  \\
6   & \multicolumn{2}{c|}{$(2,2)$}  &$(2,1)$ & &
 $(1,1)$  &$(2,0)$ &  & $(1,0)$  &  \\
5   &  \multicolumn{2}{c|}{$N=8$ }  &$N=6$   & & \multicolumn{2}{c|}{
$N=4$ }&  & $N=2$  &  \\
4  &  \multicolumn{2}{c|}{$N=8$ }  & $N=6$ & $N=5$ & \multicolumn{2}{c|}{
  $N=4$  }  &$N=3$  &  $N=2$   &  $N=1$  \\
\hline   & \multicolumn{4}{c|}{SUGRA}  &
 \multicolumn{2}{c|}{SUGRA/SUSY} & SUGRA & \multicolumn{2}{c|}{SUGRA/SUSY}  \\
\hline
\end{tabular}
\end{center}
\end{table}
%%%%%%%%%%%%%%%%%%%%%%%%%%%%%%%%%%%%%%%%%%%%%%%%%%%%%%%%%%%%%%%%

Nearly all these theories have scalar fields, which are maps from
spacetime to a `target space'. These target spaces have interesting
geometrical properties. These geometries in the case of more than 8 real
supercharges are shown in table \ref{tbl:geometriesPlus8}.
\begin{table}[htbp]
  \caption{{\it Scalar geometries in theories with more than 8
  supersymmetries (and dimension $\geq 4$).} The theories are ordered as
  in table~\ref{tbl:mapsusy}. For more than 16 supersymmetries, there is a
  unique supergravity (up to gaugings irrelevant to the geometry),
  while for 16 and 12 supersymmetries there is a number $n$ indicating
  the number of vector multiplets that are included.}\label{tbl:geometriesPlus8}
\begin{center}
  $\begin{array}{|@{\hspace{2pt}}c| % dimensie
  c| @{\hspace{2pt}}c| %einde 32
c|@{\hspace{2pt}}c| c| @{\hspace{2pt}}c| *{1}{c|}
   }
\hline
 D & \multicolumn{2}{c|}{32} & 24 &20  & \multicolumn{2}{c|}{16} & 12      \\
\hline
 10 & \SO(1,1)&\frac{\SU(1,1)}{\U(1)}& &&&&\\
 9  & \multicolumn{2}{c|}{\frac{\Sl(2)}{\SO(2)} \otimes
  \SO(1,1)}&
&& \frac{\SO(1,n)}{\SO(n)} \otimes  \SO(1,1)& &    \\[5mm]
8  & \multicolumn{2}{c|}{\frac{\Sl(3)}{\SO(3)}\otimes
\frac{\Sl(2)}{\SO(2)} }& & & \frac{\SO(2,n)}{\SO(2)\times
  \SO(n)}\otimes  \SO(1,1)&  &      \\[5mm]
7  & \multicolumn{2}{c|}{\frac{\Sl(5)}{\SO(5)} }& & &
\frac{\SO(3,n)}{\SO(3)\times
  \SO(n)}\otimes  \SO(1,1)& &      \\[5mm]
6  & \multicolumn{2}{c|}{\frac{\SO(5,5)}{\SO(5)\times \SO(5)}}
&\frac{\SO(5,1)}{\SO(5)} & & \frac{\SO(4,n)}{\SO(n)\times \SO(4)}\otimes
\SO(1,1)
 & \frac{\SO(5,n)}{\SO(n)\times \SO(5)} &     \\[5mm]
5  & \multicolumn{2}{c|}{\frac{\mathrm{E}_{6,6}}{\USp(8)}} &
\frac{\SU^*(6)}{\USp(6)} &  &
\multicolumn{2}{c|}{\frac{\SO(5,n)}{\SO(5)\times
  \SO(n)}\otimes  \SO(1,1)}   &      \\[4mm]
4  & \multicolumn{2}{c|}{\frac{\mathrm{E}_{7,7}}{\SU(8)}} &
\frac{\SO^*(12)}{\U(6)} & \frac{\SU(1,5)}{\U(5)} &
\multicolumn{2}{c|}{\frac{\SU(1,1)}{\U(1)}\times \frac{\SO(6,n)}
{\SO(6)\times \SO(n)}}  & \frac{\SU(3,n)}{\U(3)\times \SU(n)}    \\
\hline
\end{array}$
\end{center}
\end{table}
One notices that they are all symmetric spaces. On the other hand, the
theories with 4 real supersymmetries, which are the $N=1$ theories in 4
dimensions, lead to general K{\"a}hler manifolds. In both these cases, we
thus obtain geometric structures that are well known. The geometrically
interesting case are the supersymmetric theories with 8 real
supercharges. The type of geometries depends on the occurrence of
different representations (multiplets) of the supersymmetry algebra. For
our purpose here, we can restrict our attention to vector multiplets and
hypermultiplets. Vector multiplets in 6 dimensions do not have scalars
and thus no associated target-space geometry. Vector multiplets in 5
dimensions have real scalars that parametrize a geometry that is denoted
by `very special real geometry'. Those in 4 dimensions have complex
scalars that parametrize a restricted class of K{\"a}hler geometries denoted
by `special K{\"a}hler geometry'. Furthermore, one can have hypermultiplets
in dimensions $D=6,5,4$ and 3. In dimensions $D=6,5,4$ their scalars
parametrize a quaternionic-K{\"a}hler manifold, while in $D=3$ one can have a
direct product of 2 quaternionic-K{\"a}hler manifolds, which is essentially
due to the fact that the `would be' vector multiplets in $D=3$ occur as
independent hypermultiplets in that case.

All these geometries for the case of 8 supersymmetries exist in two
different versions: one which applies to rigid supersymmetry and one to
supergravity. A schematic overview of these possibilities, for $D=4$ and
$D=5$, is given in table \ref{tbl:SpecGeom}.
%%%%%%%%%%%%%%%%%%%%%%%%%%%%%%%%%%%%
\begin{table}[htbp]
  \caption{\it Geometries from supersymmetric theories with 8 real
  supercharges with vector multiplets and hypermultiplets.   }\label{tbl:SpecGeom}
\begin{tabular}{|c|ccc|}
\hline
   & $D=5$ vector multiplets & $D=4$ vector multiplets & hypermultiplets \\
\hline
 rigid & affine  & affine &  \\
 (affine) & very special real & special K{\"a}hler & hyper-K{\"a}hler  \\ \hline
 local & (projective) & (projective) &   \\
 (projective) & very special real & special K{\"a}hler  & quaternionic-K{\"a}hler \\
\hline
\end{tabular}
\end{table}
%%%%%%%%%%%%%%%%%%%%%%%%%%%%%%%%%%%%%
As indicated there, the geometries appearing in supergravity can be
considered as projective versions of the `affine' ones in rigid
supersymmetry. When one uses the terminology `very special real
manifolds' or `special K{\"a}hler manifolds', one usually refers to the
versions in supergravity. Very special real manifolds were first found
in~\cite{Gunaydin:1984bi} and connected to special geometry
in~\cite{deWit:1992cr}, where they got their name. Special K{\"a}hler
geometry was found in \cite{deWit:1984pk}, and denoted as such in
\cite{Strominger:1990pd}. A coordinate-independent formulation was found
in \cite{Castellani:1990zd,D'Auria:1991fj}. The version in rigid
supersymmetry was first investigated
in~\cite{Sierra:1983cc,Gates:1984py}. This is called e.g. `rigid special
K{\"a}hler' or was appropriately called `affine special K{\"a}hler' in
\cite{Alekseevsky:1999ts} when the supergravity version is called
projective special K{\"a}hler\footnote{A mathematical definition of very
special K{\"a}hler geometry can be found in \cite{Alekseevsky:2001if}.
Definitions of special K{\"a}hler manifold independent of supergravity were
given in~\cite{Craps:1997gp}, and a review appeared
in~\cite{VanProeyen:1999ya}. Other mathematical definitions of special
K{\"a}hler geometries have been given
in~\cite{Freed:1997dp,Alekseevsky:1999ts}.}. In that sense if
hyper-K{\"a}hler geometry and quaternionic-K{\"a}hler geometry would not have got
these names before, appropriate names would be affine, respectively
projective quaternionic-K{\"a}hler geometry. This is in fact what we will
clarify in this paper, using methods that connect also the manifolds of
the other two columns of table~\ref{tbl:SpecGeom}.

Complex structures are endomorphisms $J$ on the tangent space that square
to $-\unity $, and are 1-integrable. A hypercomplex structure has 3 such
operations with $J^1J^2=J^3$, which we collectively denote by $\vec{J}$.
Hermitian metrics $g$ obey $g(JX,JY)=g(X,Y)$. This leads to the following
characterization of manifolds that we mentioned here, which all have an
hermitian metric:
\begin{eqnarray}
 \mbox{K{\"a}hler manifolds} & : & \mbox{complex structure with }\nabla J=0, \nonumber\\
 \mbox{hyper-K{\"a}hler man.} & : & \mbox{hypercomplex structure with
 }\nabla\vec{J}=0, \nonumber\\
 \mbox{quaternionic-K{\"a}hler}&:& \mbox{hypercomplex structure with }
 \nabla \vec{J}+2\vec{\omega }\times\vec{J}=0.
 \label{defKahlerquat}
\end{eqnarray}
The first two involve the Levi-Civita connection, while the latter
condition involves moreover an $\SU(2)$ connection 1-form $\vec{\omega}$.
The $\times $ symbol in this equation is the exterior product in the
3-dimensional vector space.

Up till now all the geometries were based on a manifold with a metric,
which is in the physical theory related to the existence of a Lagrangian.
The dynamical equations of motion then follow from Euler-Lagrange
equations. However, in supersymmetric theories the dynamical equations
may also be determined by the supersymmetry algebra. This leads to
theories where the dynamics is only governed by field equations rather
than by an action. To illustrate the difference, consider the action
\begin{equation}
  S=\int \rmd t\ {\cal L}=\int \rmd t\ g_{ij}(\phi )\frac{\rmd \phi ^i}{\rmd t}\frac{\rmd \phi ^j}{\rmd
  t}.
 \label{Sphi}
\end{equation}
The Euler-Lagrange equation becomes the geodesic equation
\begin{equation}
  \frac{\rmd^2\phi ^i}{\rmd t^2}+\Gamma ^i_{jk}(\phi )\frac{\rmd \phi ^j}{\rmd t}\frac{\rmd \phi ^k}{\rmd
  t}=0.
 \label{ELphi}
\end{equation}
Note that while (\ref{Sphi}) involves a metric, the geodesic equation
involves only an affine connection. In this case, the affine connection
is the Levi-Civita connection, but in general one could consider
(\ref{ELphi}) with another affine connection. In the applications that we
have in mind, the closure of the supersymmetry algebra leads directly to
equations similar to (\ref{ELphi}), which do not necessarily involve a
metric. In the above equation only a torsionless connection occurs. We
will show below that in the supersymmetric theories that we consider we
have to require the affine connection to be torsionless.

\section{Quaternionic-like manifolds}
 \label{ss:quatlike}
\subsection{Affine connections}
 \label{ss:affineconnections}
We will now repeat some properties of the family of quaternionic-like
manifolds. As we mentioned above, hypercomplex structures are defined by
3 endomorphisms denoted as $H=\{\vec{J}\}=\{J^1,J^2,J^3\}$. A
quaternionic structure is the linear space $Q=\{\vec{a}\cdot
\vec{J}|\vec{a}\in\mathbb{R}^3\}$. A hermitian bilinear form is a form
$F(X,Y)$ with $F(JX,JY)=F(X,Y)$. A non-singular hermitian bilinear form
is a `good metric'. As such we can define a quartet of quaternionic-like
manifolds in table \ref{tbl:quatlikeMan}. The table is essentially taken
over from~\cite{AM1996}, where even more distinction has been made
between various cases. The quaternionic manifolds represent the generic
case that includes all the others as special cases.
\begin{table}[ht]
  \caption{\it Quaternionic-like manifolds of real dimension $4r$, and their holonomy groups.}
  \label{tbl:quatlikeMan}
\begin{center}
  \begin{tabular}{|c||c|c||l|}\hline
      & no good  metric & with a good metric& \\ \hline\hline
    no $\SU(2)$ & \textit{\textbf{hypercomplex}} & \textit{\textbf{hyper-K{\"a}hler}}& rigid  \\
    connection & $\Aut(H)=$& $\Aut(H,g)=$&supersymmetry \\
& $\Gl(r,\mathbb{H})$ & $\USp(2r)$& \\ \hline
    non-zero $\SU(2)$ & \textit{\textbf{quaternionic}} & \textit{\textbf{quaternionic-K{\"a}hler}}& \\
    connection & $\Aut(Q)=$& $\Aut(Q,g)=$ & supergravity \\
   &$ \SU(2)\cdot \Gl(r,\mathbb{H})$  & $\SU(2)\cdot\USp(2r)$ &  \\
\hline\hline
 & field equations & action& \\ \hline
  \end{tabular}
\end{center}
\end{table}

A hypercomplex manifold is equipped with an affine connection such that
\begin{equation}
  \nabla \vec{J}=0.
 \label{hypercNabla}
\end{equation}
Given the hypercomplex structure, this connection is unique and is in
general the sum of a so-called `Obata connection' and the Nijenhuis
tensor. The last part is the torsion, and as we mentioned that we are
interested in torsionless connections, it should vanish.

For a quaternionic manifold we only need
\begin{equation}
  \nabla Q\subset Q, \qquad \mbox{i.e.}\qquad \nabla \vec{J}+2\vec{\omega
  }\times\vec{J}=0.
 \label{quatNabla}
\end{equation}
In this case, the connection is not unique, even not for torsionless
connections to which we will restrict here. Indeed, for a solution of
(\ref{quatNabla}) and any 1-form $\xi $ we can construct other solutions
of (\ref{quatNabla}) as
\begin{equation}
 \nabla '=\nabla + S^\xi ,\qquad \vec{\omega }'=\vec{\omega
 }+\vec{J}^*\xi,
 \label{xitransf}
\end{equation}
with
\begin{equation}
 S^\xi _X Y= \xi (X) Y+\xi (Y) X -\xi (\vec{J}X)\cdot \vec{J}Y-\xi
(\vec{J}Y)\cdot
 \vec{J}X.
\end{equation}
Of course, when we are discussing hyper-K{\"a}hler or quaternionic-K{\"a}hler
manifolds, the affine connections should be the Levi-Civita connections.
The condition (\ref{quatNabla}) is a weaker condition than
(\ref{hypercNabla}). The condition that this can be solved is that the
Nijenhuis tensor is of the form
\begin{equation}
  N_{XY}{}^Z=-\vec J_{[X}{}^Z \cdot \vec \omega^{\rm Op} _{Y]}=
  -\ft12\vec J_{X}{}^Z \cdot \vec \omega^{\rm Op} _{Y}+\ft12\vec J_{Y}{}^Z \cdot \vec \omega^{\rm Op} _{X},
 \label{NijenhQuat}
\end{equation}
where $\vec \omega^{\rm Op}$ is at this point an arbitrary triplet of
1-forms. $X,Y,Z$ in this equation are indices labelling the $4r$
coordinates $q^X$ of the manifold. E.g. $\vec{\omega }=\vec{\omega
}_X\rmd q^X$. We warn the reader that $X$ has been used also as
indication of a vector field, and it should be clear from the context
what is meant. $\vec \omega^{\rm Op}$ is an $\SU(2)$ connection such that
the Nijenhuis condition (\ref{NijenhQuat}) guarantees the existence of
corresponding affine connection coefficients
\begin{equation}
  \Gamma ^{\rm Op}{}_{XY}{}^Z\equiv \Gamma^{\rm
Ob}{}_{XY}{}^Z-\ft12\vec J_X{}^Z\cdot  \vec \omega^{\rm Op} _Y-\ft12\vec
J_Y{}^Z\cdot  \vec \omega^{\rm Op} _X,
 \label{GammaOp}
\end{equation}
where the first term involves the Obata connection coefficients. This
solves (\ref{quatNabla}) and is called the Oproiu connection. The
connections for the different quaternionic manifolds are schematically
shown in table \ref{tbl:ConnQuatlikeMan}.
\begin{table}[t]
  \caption{\it The affine connections in quaternionic-like manifolds}
\label{tbl:ConnQuatlikeMan}
\begin{center}
  \begin{tabular}{||c|c||}\hline\hline
   \textit{\textbf{hypercomplex}} & \textit{\textbf{hyper-K{\"a}hler}}  \\
    Obata connection & Obata connection \\
                     & = Levi-Civita connection  \\ \hline
    \textit{\textbf{quaternionic}} & \textit{\textbf{quaternionic-K{\"a}hler}} \\
    Oproiu connection or  & Levi-Civita connection = \\
    other related by $\xi _X$ transformation & connection related to
    Oproiu \\
    & by a particular choice of $\xi _X$ \\
\hline\hline
  \end{tabular}
\end{center}
\end{table}
Note that the freedom of choice for $\xi $ disappears in
quaternionic-K{\"a}hler manifolds, where the $\xi$ is determined by the
requirement that the affine connection should coincide with the
Levi-Civita connection.

\subsection{Supersymmetry}

We will illustrate here how the hypercomplex geometry arises in
supersymmetric models, and why this needs torsionless connections. This
statement holds\footnote{Hyper-K{\"a}hler manifolds with torsion do appear in
2-dimensional supersymmetric theories \cite{Gates:1984nk,Howe:1987qv}.}
for supersymmetric theories in spacetime dimensions $D=3,4,5$ or 6. The
essence of this part is not dependent on whether we consider any of these
dimensions, but the notation is simplest for $D=5$ or 6, to which we will
restrict ourselves for convenience. The hypermultiplets consist of fields
$q^X(x)$ that are maps from spacetime with coordinates $x^a$ to the
quaternionic space with coordinates $q^X$, with $X=1,\ldots ,4r$, and
fermionic partners, which are spinors of spacetime $\zeta ^A$(x), where
$A=1,\ldots ,2r$, and spinor indices are suppressed.  They belong to the
irreducible spinor modules for $\SO(4,1)$, resp. $\SO(5,1)$. They are
`symplectic Majorana spinors' for $D=5$ or `symplectic Majorana-Weyl' for
$D=6$, using the symplectic matrix $\rho^{AB}$, which satisfies
\begin{equation}
  \rho ^{AB}=-\rho ^{BA},\qquad \rho _{AB}=(\rho ^{AB})^*,\qquad
  \rho ^{AB}\rho _{CB}=\delta ^A_C.
 \label{proprho}
\end{equation}
We suppress below the dependence of the fields $q^X$ and $\zeta ^A$ on
spacetime coordinates.

Supersymmetry is defined also by a symplectic spinor generator $Q_i$ (for
$i=1,2$, related to the notation $N=2$ in table \ref{tbl:mapsusy}), where
the symplectic matrix is $\varepsilon^{ij}=\varepsilon _{ij}$, with
$\varepsilon ^{12}=-\varepsilon ^{21}=1$, satisfying the same relations
as in (\ref{proprho}). The supersymmetry transformations are denoted as
\begin{equation}
  \delta (\epsilon )=\bar \epsilon ^iQ_i,
 \label{susytransf}
\end{equation}
where $\epsilon^i$ are the parameters of the supersymmetry
transformations and $\bar \epsilon ^i$ is the Majorana conjugate.

The supersymmetry algebra is
\begin{equation}
  \left[ \delta (\epsilon _1),\delta (\epsilon _2)\right] = \ft12\bar
  \epsilon_2 ^i \gamma ^a \epsilon _1^j \varepsilon _{ji}P_a,   \qquad
  a=1,\ldots D.
 \label{susyalgebra}
\end{equation}
where $\gamma^a$ are the matrices used for the Clifford algebra, and
satisfy $\gamma ^a\gamma ^b+\gamma ^b\gamma ^a=2\eta ^{ab}$ where $\eta
=\diag(-1,+1,\ldots ,+1)$. The generators $P_a$ are translations that act
on the $q^X$ as $P_aq^X=\partial _aq^X$, where $\partial_a$ are the
derivatives with respect to the spacetime coordinates $x^a$.

The supersymmetry transformations of the $q^X$ take the general form
\begin{equation}
  \delta (\epsilon )q^X\equiv  \bar \epsilon ^iQ_iq^X= -\rmi\bar \epsilon ^i
  \zeta ^A f_{iA}^X(q),
 \label{susyq}
\end{equation}
This implies for consistency reality conditions on the coefficients
functions $f_{iA}^X(q)$:
\begin{equation}
  (f_{iA}^X)^*= \varepsilon ^{ij}\rho ^{AB}f_{jB}^X,
 \label{realityf}
\end{equation}

To generate a translation in the commutator on $q^X$ according to
(\ref{susyalgebra}), the supersymmetry transformation of the spinor must
contain a term of the form
\begin{equation}
  \delta (\epsilon )\zeta ^A=\ft12\rmi \gamma ^a \epsilon _i\left[ \partial _a q^X\right]
  f^{iA}_X(q)+\ldots ,
 \label{susyzeta}
\end{equation}
where
\begin{equation}
  f_X^{iA}f_{iA}^Y=\delta _X^Y,
 \label{finvers}
\end{equation}
i.e.\ the $f^{iA}_X$ and $f^X_{iA}$ are each others inverse as $4r\times
4r$ matrices. Note that we also suppress the dependence on the
coordinates $q^X$. However, when calculating the commutator on $q^X$ we
have to take into account the $q$-dependence of $f^X_{iA}$, leading to a
term
\begin{equation}
  \left[ \delta (\epsilon _1),\delta (\epsilon _2)\right]q^X=\ldots -\rmi\bar \epsilon_2 ^i
  \zeta ^A \left[ \partial _Yf_{iA}^X \right] \left[ \delta (\epsilon _1)q^Y\right]
   - (\epsilon _1\leftrightarrow\epsilon _2).
 \label{calccomm}
\end{equation}
To remove this term, we can modify the transformation of $\zeta ^A$ by a
term proportional to another $\zeta ^B$ and a supersymmetry
transformation of $q^X$, i.e.\ we complete (\ref{susyzeta}) to a form
\begin{equation}
  \delta (\epsilon )\zeta ^A= \ft12\rmi\gamma ^a \epsilon _i\left[ \partial _a q^X\right]
  f^{iA}_X-\zeta ^B\omega _{XB}{}^A(q)\left[ \delta (\epsilon )q^X\right]   ,
 \label{delsusyzeta}
\end{equation}
where $\omega _{XB}{}^A(q)$ has to be determined. Using (\ref{susyq}) in
(\ref{calccomm}) and adding the contribution of the last term of
(\ref{delsusyzeta}), leads to
\begin{eqnarray}
  \left[ \delta (\epsilon _1),\delta (\epsilon _2)\right]q^X&=\ft14\bar
  \epsilon_2 ^i \gamma ^a \epsilon _1^j \varepsilon _{ji}\partial _aq^X &+ f^{iA}_Z
  \left[ \delta (\epsilon _2)q^Z\right]
  \partial _Yf_{iA}^X \left[ \delta (\epsilon _1)q^Y\right]\label{calccomm2}\\
&&-\left[ \delta (\epsilon _2)q^Z\right]f_Z^{iB} f_{iA}^X\omega
_{YB}{}^A\left[ \delta (\epsilon_1 )q^Y\right]
   - (\epsilon _1\leftrightarrow\epsilon _2).\nonumber
\end{eqnarray}
In order that the latter terms do not contribute to the commutator, they
should add to a symmetric expression in $(YZ)$ for any $X$, which we
denote as $\Gamma ^X_{YZ}$:
\begin{equation}
  f^{iA}_Z  \partial _Yf_{iA}^X -f_Z^{iB} f_{iA}^X\omega
_{YB}{}^A=-\Gamma _{YZ}^X.
 \label{reqvanishcomm}
\end{equation}
This equation is equivalent to the requirement
\begin{equation}
  \partial _Yf_{iA}^X-\omega _{YA}{}^Bf_{iB}^X +f_{iA}^Z\Gamma _{YZ}^X=0.
 \label{covconstf}
\end{equation}
This is the condition of covariant constancy of a vielbein $f_{iA}^X$ in
the quaternionic manifold with a torsionless connection $\Gamma
_{YZ}^X=\Gamma _{ZY}^X$ and $\omega_{XA}{}^B$ is a $\Gl(r,\mathbb{H})$
connection (written as $2r\times 2r$ complex matrices with reality
properties determined by the consistency in the transformation
(\ref{delsusyzeta})).

We want to remark that in the case of local supersymmetry, this condition
can be relaxed. Indeed, in this case the algebra can contain a
field-dependent supersymmetry transformation. This means that
(\ref{reqvanishcomm}) can include in the right-hand side a term of the
form $-f^X_{iA}f_Z^{jA}\omega _{Yj}{}^i(q)$, where $\omega _{Yj}{}^i(q)$
is arbitrary and defines an $\SU(2)$ connection. This possibility does
not apply to rigid supersymmetry as it involves a field-dependent
supersymmetry transformation.

\subsection{Vielbeins}

The supersymmetry analysis of the previous subsection leads to the result
that for rigid supersymmetry one needs a vielbein that satisfies the
integrability condition (\ref{covconstf}), including a torsionless
connection and a connection for $\Gl(r,\mathbb{H})$.  The vielbeins
determine the hypercomplex structure as
\begin{equation}
  \vec J_X{}^Y=-{\rm i}\vec \sigma _i{}^jf_X^{iA}f^Y_{jA}.
 \label{defJfromf}
\end{equation}
The factors $f_X^{iA}f^Y_{jA}$ define for any $X$, $Y$, a $2\times 2$
matrix with trace $\delta _X^Y$. The $\vec{\sigma }$ are the traceless
Hermitian Pauli matrices and project the other 3 components of the
$2\times 2$ matrix. These expressions automatically satisfy the relations
for a hypercomplex structure. Moreover, the covariant constancy of the
vierbein according to (\ref{covconstf}) implies the covariant constancy
(\ref{hypercNabla}) of the hypercomplex structure.

The condition (\ref{covconstf}) can always be solved for the connection
$\omega _{XA}{}^B$ once we know the vielbein and the $\Gamma _{YZ}^X$. As
mentioned above, the torsionless connection can be found if the Nijenhuis
tensor vanishes, and is in that case uniquely determined.

For quaternionic manifolds, the vielbeins should
satisfy\footnote{\label{ijtoalpha}One can make the transition from
doublet to vector notation by using the sigma matrices, ${\omega}_{X
i}{}^j = \rmi \vec \sigma _i{}^j\cdot \vec \omega_X$, and similarly $\vec
\omega_X = -\frac 12\rmi\vec \sigma _i{}^j {\omega}_{Xj}{}^i$. This
transition between doublet and triplet notation is valid for any triplet
object as e.g.\ the complex structures.}
\begin{equation}
  \partial _X f_Y^{iA}-\Gamma _{XY}^Zf_Z^{iA}+f_Y^{jA}\omega _{Xj}{}^i +f_Y^{iB}\omega
_{XB}{}^A=0.
 \label{covconstfSU2}
\end{equation}
Also here, this equation determines $\omega _{XA}{}^B$ once we know the
other connections. But now the latter are not uniquely defined by
(\ref{quatNabla}), but allow the $\xi$-transformations (\ref{xitransf}).

\subsection{Curvature decompositions}
 \label{ss:curvdecomp}

The relation (\ref{covconstfSU2}) has as integrability condition a
relation between curvatures:
\begin{eqnarray}
  \Purple{R_{XYW}{}^Z}&=\OliveGreen{R^{\SU(2)}{}_{XYW}{}^Z}&+\,\Red{R^{\Gl(r,\mathbb{H})}{}_{XYW}{}^Z}
  \nonumber\\
&=-\vec J_W{}^Z\cdot\OliveGreen{\vec {\cal R} _{XY}}& +\,
L_W{}^Z{}_A{}^B\,\Red{{\cal R}_{XYB}{}^A }. \label{curvatrel}
\end{eqnarray}
The left-hand side is the curvature defined by the affine connection
$\Gamma _{YZ}^X$, while $\vec {\cal R} _{XY}$ is the curvature determined
by the $\SU(2)$ connection $\vec{\omega }_X$ and ${\cal R}_{XYB}{}^A $ is
determined by $\omega _{XB}{}^A$. The object $L_W{}^Z{}_A{}^B$ is defined
similar to the complex structures (\ref{defJfromf}), but with contraction
over indices $i$ rather than $A$:
\begin{equation}
  L_W{}^Z{}_A{}^B\equiv f^Z_{iA}f_W^{iB}.
 \label{defL}
\end{equation}
The relation (\ref{curvatrel}) holds for general quaternionic manifolds.
In the case of hypercomplex or hyper-K{\"a}hler manifolds, the $\SU(2)$ term
is absent. The Ricci tensor ${\it Ric}$ determined by the curvature
$R_{XYW}{}^Z$ can in general have an antisymmetric part if the trace of
the $\Gl(r,\mathbb{H})$ curvature is non-vanishing:
\begin{equation}
  \Purple{Ric_{[XY]}}=\Purple{R_{Z[XY]}{}^Z}=-\Red{{\cal R}_{XY}^{\mathbb{R}}}\equiv -
   \Red{{\cal R}_{XYA}{}^A}.
 \label{RicASfrom}
\end{equation}
As the Ricci tensor associated to a Levi-Civita connection is symmetric,
this has to vanish in the case of hyper-K{\"a}hler or quaternionic-K{\"a}hler
manifolds.

An unnatural feature of the splitting (\ref{curvatrel}) is that the
individual terms do not satisfy the first Bianchi identity. An
alternative splitting is
\begin{equation}
  R_{XYZ}{}^W=\Blue{R^{\rm Ric}{}_{XYZ}{}^W} - \frac 12 f^{Ai}_X
\varepsilon_{ij}f^{jB}_Yf_W^{kC}f^Z_{kD}
  \Magenta{\mathcal{W}_{ABC}{}^D}.
 \label{Rsecondsplit}
\end{equation}
Here both curvature tensors do satisfy the first Bianchi identity. The
first term is called the Ricci part because it is determined only by the
Ricci tensor. The second part is called the Weyl part. The Ricci tensor
of the Weyl part vanishes. It is determined by a tensor
$\mathcal{W}_{ABC}{}^D$ that is symmetric in its lower indices and
traceless. The proofs of these statements are reviewed in Appendix~B of
\cite{Bergshoeff:2002qk}.

We can summarize these curvature decompositions in the following scheme:
\begin{equation}
  \begin{array}{rccccc}
  R_{XYZ}{}^W= & \big( R^{\rm Ric}_{\rm symm}  & + & R^{\rm Ric}_{\rm antis}  & + &
  R^{({\rm W})}\big)_{XYZ}{}^W
  \\
  &\multicolumn{3}{c}{\phantom{.}\hspace{5mm}\leavevmode
\ifpdf
  \includegraphics{pijlen.pdf}
 \else
  \epsfxsize=85mm
 \epsfbox{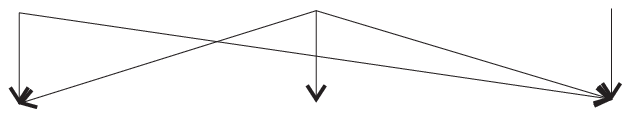}
 \fi
 \hspace{-18mm}} \\
  = & \big( R^{\SU(2)} & + & R^{\mathbb{R}} & + &
R^{\Sl(r,\mathbb{H})}\big)_{XYZ}{}^W.
\end{array}
 \label{RelatDecompR}
\end{equation}
In this decomposition we made two further splits. The Ricci curvature has
been separated in a part determined by the symmetric part of the Ricci
tensor, and a part determined by its antisymmetric part. On the other
hand, the $\Gl(r,\mathbb{H})$ has been split into $\mathbb{R}\times
\Sl(r,\mathbb{H})$.

The terms in the second line depend only on specific terms of the first
line as indicated by the arrows. This is the general scheme and thus
applicable for quaternionic manifolds, which is the general case. For
hypercomplex manifolds, we mentioned already that there is no $\SU(2)$
part, and in the upper decomposition there is no symmetric part. In this
case, the full curvature can be written as
\begin{equation}
   R_{XYZ}{}^W= - \frac 12 f^{Ai}_X
\varepsilon_{ij}f^{jB}_Yf_W^{kC}f^Z_{kD}
  \Magenta{W_{ABC}{}^D},
 \label{curvhypercomplex}
\end{equation}
where the trace  $W_{ABC}{}^C$ determines the antisymmetric Ricci tensor,
and the traceless part of $W_{ABC}{}^D$ is ${\cal W}_{ABC}{}^D$ that
appears in (\ref{Rsecondsplit}).

On the other hand, in quaternionic-K{\"a}hler manifolds there can be no
antisymmetric Ricci part and no $\mathbb{R}$-curvature. For hyper-K{\"a}hler
manifolds, both restrictions apply, and one has only the right-most terms
in both decompositions.

The $\xi$-transformation can always be used to choose connections such
that the $\mathbb{R}$-curvature vanishes.

Finally, we want to make a remark about the curvature of
quaternionic-K{\"a}hler manifolds. These are Einstein manifolds with
\begin{equation}
  Ric_{XY}=\nu (r+2) g_{XY},\qquad \OliveGreen{\vec{\cal R}_{XY}{}} =\ft12\nu
\vec J_{XY},
 \label{quatKahlerEinstein}
\end{equation}
where $\nu $ is an arbitrary real number. In supergravity this number is
related to Newton's gravitational constant:
\begin{equation}
  \nu =-\kappa ^2=-8\pi G_N.
 \label{relationnuNewton}
\end{equation}
Therefore only negative values of $\nu $ appear in supergravity, and the
scalar curvature is negative. This implies that the manifold is
non-compact (if there is at least one isometry).

\section{Conformal symmetry and the map}
 \label{ss:map}

Constructing supergravity theories is more complicated than constructing
rigid supersymmetric theories. There exists a method to construct
supergravity theories that starts from the rigid theories. These rigid
theories should be invariant under superconformal transformations. Then
one can gauge the superconformal group, and afterwards break explicitly
the symmetries that are extra with respect to the super-Poincar{\'e} group.
Indeed, the final goal is only to obtain theories that are invariant
under this latter group. However, the construction via the superconformal
group simplifies the calculations due to the larger amount of symmetry.

The gauge-fixing procedure consists in choosing a parametrization such
that for every extra symmetry there is a unique field that transforms
under it. Whenever this happens, the statement of symmetry is just that
this field is irrelevant. The remaining fields are then the physical
fields, and they do not feel the extra symmetries. The main extra
symmetry is the dilatation. We will choose one field that describes the
scale, which is the field that is fixed in the procedure described above.
The superconformal group includes for our case also an $\SU(2)$ group
which we will use to eliminate a further 3 fields. Mathematically this
means that we consider first a projective version of the theory that we
want to describe.

A conformal symmetry amounts to the presence of a vector $k$ such that
the Lie derivative of the metric is proportional to the metric: ${\cal
L}_k g=w\,g$ for a constant (positive) $w$. Such a vector is called a
homothetic Killing vector. Special conformal transformations need an
extra condition, namely that the one-form $g(k,\cdot )$ is closed. The
combination of these two conditions can be written as
\begin{equation}
  \nabla _X k =\ft12 w X.
 \label{closedhomoth}
\end{equation}
A vector that satisfies this condition is `a closed homothetic Killing
vector'. We use the normalization $w=3$. It is important to notice that
the condition (\ref{closedhomoth}) is independent of a metric. Therefore
we can use the concept of closed homothetic Killing vectors for manifolds
without a good metric, despite the fact that homothetic Killing vectors
are only defined with respect to a metric.

The presence of the hypercomplex structure implies that the vectors
\begin{equation}
  \vec{k}=\ft13\vec{J}\,k
 \label{defveck}
\end{equation}
generate an $\SU(2)$, which is the subgroup of the superconformal group
mentioned above.

The general strategy is thus to start with a hypercomplex manifold that
has a closed homothetic Killing vector. We will denote such a manifold as
a `conformal hypercomplex manifold'. Assume that the dimension of this
manifold is $4(n_H+1)$. Then we will isolate 4 directions in this
manifold that transform under the dilatation and the $\SU(2)$
transformations, and the orthogonal $4n_H$ dimensional manifold that is
invariant. This invariant submanifold then inherits the property of being
quaternionic. Some of these steps have already been performed in
\cite{deWit:1999fp}. Recently, we \cite{Bergshoeff:2004nf} clarified in
this way the general structure of the map, especially showing its one to
one character. This we will review below.

\subsection{Conformal hypercomplex manifolds}

Following the previous ideas, it is appropriate to choose coordinates
adapted to the conformal structure. As such we define a first coordinate
that we will denote by $z^0$ such that the closed homothetic Killing
vector points in this direction. Denoting the coordinates of this
hypercomplex manifold by $q^{\hat{X}}$, we choose
\begin{equation}
  k^{\hat{X}}=\delta^{\hat{X}}_0 k^0= 3\delta^{\hat{X}}_0 z^0.
 \label{normalz0}
\end{equation}
The factor 3 is purely a matter of normalization. Then we choose 3 more
coordinates such that the vectors $\vec{k}$ in (\ref{defveck}) only point
in these three directions. We denote these directions with an index
$\alpha =1,2,3$, and thus the vectors $\vec{k}$ have only nonzero
components $\vec{k}^\alpha $. All other components of $q^{\hat{X}}$ are
denoted by $q^X$:
\begin{equation}
  q^{\hat X}=\left\{ \Red{z^0, z^\alpha} , \Blue{q^X}\right\}.
 \label{allcomponents}
\end{equation}
The strategy is to fix $z^0$ by a gauge choice for dilatations, and
$z^\alpha $ by fixing the $\SU(2)$ symmetries. The gauge-fixed manifold
thus contains only the directions $q^X$, and will be the quaternionic
manifold.

This choice of coordinates and the hypercomplex algebra imply that the
hypercomplex structures decompose as
\begin{equation}
  \begin{array}{lll}
    \widehat{\vec J}_0{}^0=0 \,,\quad & \widehat{\vec J}_\alpha{}^0=-\Red{z^0\vec m_\alpha}
      \,, & \widehat{\vec J}_X{}^0=\Red{z^0} \OliveGreen{\vec A _X}  \,,\\
    \widehat{\vec J}_0{}^\beta=\frac{1}{\Red{z^0}}\Red{\vec k^\beta }  \,,&
    \widehat{\vec J}_\alpha{}^\beta=
    \Red{\vec k^\beta} \times \Red{\vec m_\alpha }
    \,,\quad&
    \widehat{\vec J}_X{}^\beta= \OliveGreen{\vec A _X}   \times \Red{\vec k^\beta}
                               +\Blue{\vec J_X{}^Z}\left(\OliveGreen{\vec A _Z}\cdot \Red{\vec k^\beta}
                                \right)   \,,    \\
    \widehat{\vec J}_0{}^Y=0 \,,& \widehat{\vec J}_\alpha{}^Y=0 \,,& \widehat{\vec J}_X{}^Y=\Blue{\vec
    J_X{}^Y}.
  \end{array}
 \label{hypercomplexconf}
\end{equation}
In this equation $\vec{m}_\alpha $ are the inverse of $\vec{k}^\alpha $
as $3\times 3$ matrices:
\begin{equation}
  \Red{\vec k^\alpha }\cdot\Red{\vec m_\beta}=\delta ^\alpha_\beta.
 \label{mkinverse}
\end{equation}
Note that (\ref{hypercomplexconf}) depends on $z^0$, $\vec{k}^\alpha$,
${\vec{A}_X}$ and $\vec{J}_X{}^Y$. We mentioned already that $z^0$ is a
scale variable and $\vec{k}^\alpha$ are the $\SU(2)$ Killing vectors.
Furthermore, there is the triplet of one-forms $\vec{A}=\vec{A}_X\rmd
q^X$, which are arbitrary up to this point and the $\vec{J}_X{}^Y$, which
satisfy the hypercomplex algebra by itself. The latter span the
quaternionic structure on the $4n_H$-dimensional submanifold.

Up to now, the matrices (\ref{hypercomplexconf}) define an almost
hypercomplex structure. In order to become an hypercomplex structure we
need to impose the vanishing of the Nijenhuis tensor
$N_{\hat{X}\hat{Y}}{}^{\hat{Z}}=0$. This leads to two further conditions:
\begin{itemize}
  \item The curvature of the triplet $\vec{A}$ is related to the complex
  structure:
\begin{equation}
  \left( 2\rmd \vec{A} -\vec{A}\times \vec{A}\right)
  (X,Y)=h(\vec{J}\,X,Y)-h(\vec{J}\,Y,X),
 \label{cond1Nijenh}
\end{equation}
where the wedge product between the one forms $\vec{A}$ is understood.
$h$ is an arbitrary symmetric bilinear form. Here the $X,Y$ denote
vectors of the $4n_H$-dimensional subspace.
  \item The subspace is quaternionic. This means that the Nijenhuis
  tensor of the $J_X{}^Y$ satisfies (\ref{NijenhQuat}), with
\begin{equation}
   \OliveGreen{ \vec \omega^{\rm Op}_X}=-\ft{1}6\left(
2\OliveGreen{\vec{A}_X}+\OliveGreen{\vec{A}_Y}\times
  \Blue{\vec J_X{}^Y}\right).
 \label{Oproiuquat}
\end{equation}
\end{itemize}
This form of the $\SU(2)$ connection gives the Oproiu connection. We know
already that we can use the freedom of $\xi $-transformations in
(\ref{xitransf}) to obtain other forms of the $\SU(2)$ connections. In
this way, we can simplify it to
\begin{equation}
 \vec \omega_X =-\ft12\vec{A}_X.
 \label{vecomegaXsimple}
\end{equation}
This choice has further advantages: the $\mathbb{R}$-curvature of the
quaternionic manifold is the same as the $\mathbb{R}$-curvature of the
hypercomplex manifold. As a consequence, in this $\xi $-choice, the
$\mathbb{R}$-curvatures of the quaternionic manifold is Hermitian,
because $\mathbb{R}$-curvatures of hypercomplex manifolds are always
Hermitian \cite{AM1996}.

\subsection{The $\hat{\xi }$-transformations}

The above analysis of the general form of the hypercomplex structures for
conformal hypercomplex manifolds leads to a new transformation
\cite{Bergshoeff:2004nf}, which is similar to the $\xi $-transformations
discussed in section \ref{ss:affineconnections}. Indeed, we can consider
changes of the triplet 1-forms $\vec{A}$ in (\ref{hypercomplexconf}),
such that (\ref{cond1Nijenh}) remains satisfied. Such changes are
determined by a 1-form on the quaternionic space, $\hat{\xi }$:
\begin{equation}
  \delta(\Red{\widehat{\xi }}) \OliveGreen{\vec{A}} = 2\vec{J}^*\Red{\widehat{\xi
  }}.
 \label{xihatonA}
\end{equation}
These induce therefore transformations on the hypercomplex structures,
preserving the hypercomplex algebra. Notice that this implies that the
$\hat{\xi }$ transformations have a different meaning than the $\xi
$-transformations. The latter do not transform the complex structures,
but only the connections.

We can use these $\hat{\xi }$ transformations to eliminate the
$\mathbb{R}$-curvature of the conformal hypercomplex manifold.

\subsection{Metric spaces}

If the conformal hypercomplex manifold allows a good metric, i.e.\ when
it is hyper-K{\"a}hler, one can show that it should be of the form
\begin{eqnarray}
\rmd \widehat s^2  & = & -\Red{\frac{(\rmd z^0)^2}{z^0}}
 +\Red{z^0}\Big\{ \Magenta{h_{XY}} \Blue{\rmd q^X \rmd q^Y}
 \nonumber\\
&&\qquad \qquad - \Red{\vec{m}_\alpha \cdot \vec{m}_\beta}  [\rmd
z^\alpha - \OliveGreen{\vec A_X}\cdot \Red{\vec k^\alpha} \Blue{\rmd
q^X}][\rmd z^\beta -\OliveGreen{ \vec A_Y}\cdot \Red{\vec k^\beta}
\Blue{\rmd q^Y}]\Big\},
 \label{metricconfhc}
\end{eqnarray}
where $h_{XY}$ are the components of the bilinear form that was
introduced in  (\ref{cond1Nijenh}).

We find furthermore  that the large space is a hyper-K{\"a}hler manifold if
and only if the submanifold is quaternionic-K{\"a}hler. This is the case if
$h_{XY}$ is Hermitian and invertible. Then
\begin{equation}
  g_{XY}=z^0h_{XY}
 \label{gquatK}
\end{equation}
is the metric on the quaternionic-K{\"a}hler manifold. We see that its scale
is determined by the choice of $z^0$. Indeed,  it determines the value
$\nu $ in (\ref{quatKahlerEinstein}):
\begin{equation}
  \nu =-\frac{1}{z^0}.
 \label{nuinz0}
\end{equation}
The $\SU(2)$ connection of the quaternionic-K{\"a}hler manifold is still
given by (\ref{vecomegaXsimple}). It can be shown that this corresponds
precisely to the $\xi$-gauge for which the affine connection coincides
with the Levi-Civita connection.

\subsection{Curvature mapping}

We can now compare the curvature decompositions discussed in
(\ref{ss:curvdecomp}) for the large and small manifolds. This leads to
the following scheme:
\begin{equation}
  \begin{array}{ccccccc}
    \widehat{R} & = &   &   & \Blue{\widehat{R}^{\rm Ric}_{\rm antis}} & + & \Magenta{\widehat{R}^{(\rm W)}} \\
            &   &   &   & \uparrow & \nearrow & \uparrow \\
     &   & h_{XY} &   & \widehat{W}_{ABC}{}^C &   & \mathcal{W}_{ABC}{}^D \\
      &   & \downarrow &   & \downarrow &   & \downarrow \\
    R & = & \Blue{R^{\rm Ric}_{\rm symm}}  & + & \Blue{R^{\rm Ric}_{\rm antis}}  & + & \Magenta{R^{(\rm W)}} \\
  \end{array}
 \label{relationsdecomp}
\end{equation}
The upper line and lower line are respectively the curvature
decompositions of the hypercomplex and quaternionic manifolds. One
remarks that $h$ determines the symmetric Ricci tensor. On the other hand
the trace $\widehat{W}_{ABC}{}^C$ determines the antisymmetric Ricci
tensors of as well the hypercomplex as the quaternionic space, while the
traceless tensor ${\cal W}_{ABC}{}^D$ contributes to the Weyl curvature
of both manifolds. For hyper-K{\"a}hler and quaternionic-K{\"a}hler manifolds,
$\widehat{W}_{ABC}{}^C=0$, and
$\mathcal{W}_{ABC}{}^D=\widehat{W}_{ABC}{}^D$, such that there are no
antisymmetric Ricci tensors.

\subsection{The picture of the map}

\begin{figure}[p]
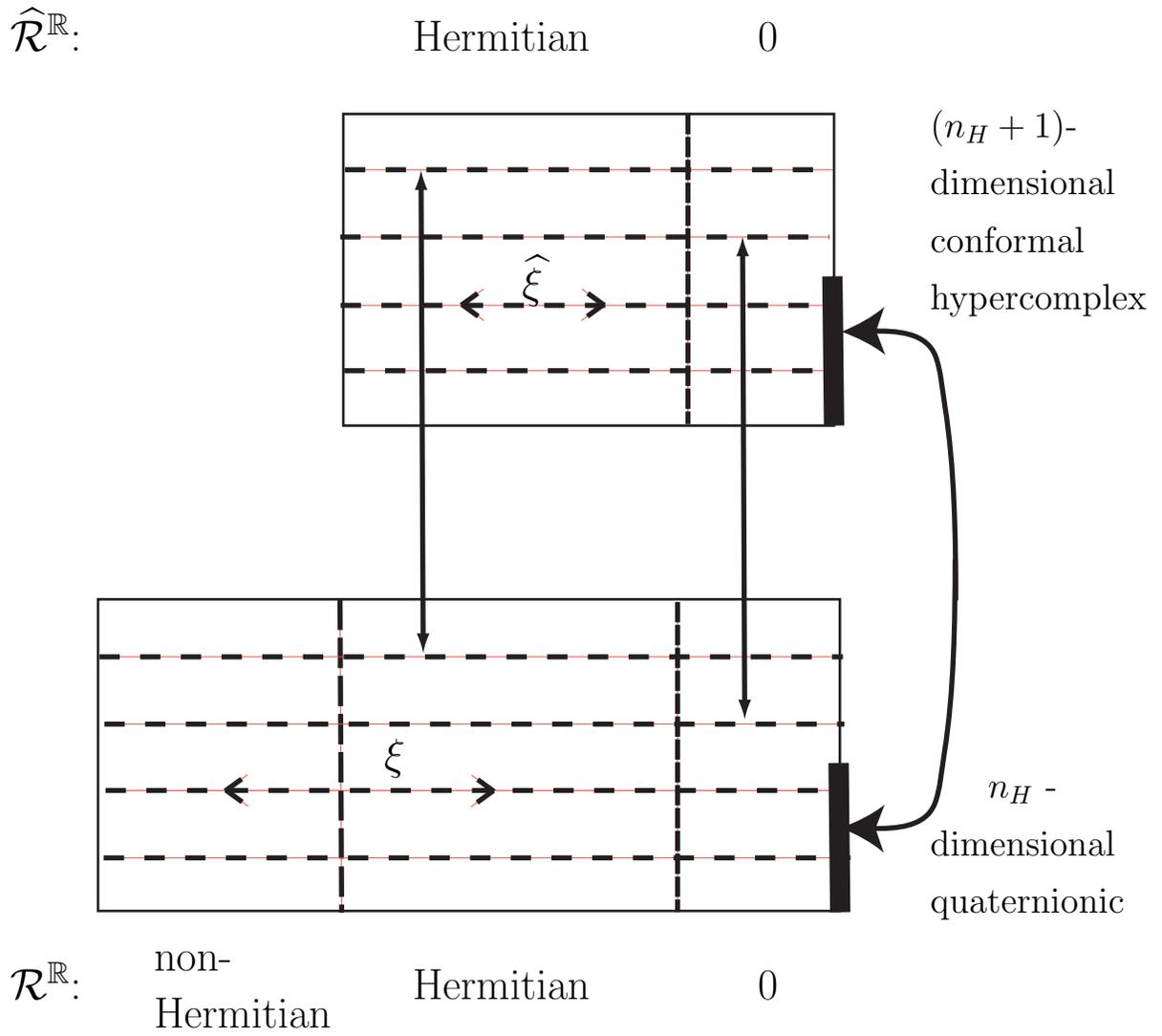

\unitlength=0.8mm \caption{\it The map schematically.}
 \label{fig:mapoverview}
\begin{center}
\begin{picture}(150,190)(0,0)
 \ifpdf
 \put(-53,-46){
  \includegraphics{maphkctoqk.pdf}}
 \else
  \put(-10,20){\leavevmode
   \epsfxsize=120mm
 \epsfbox{maphkctoqk.eps}
  } \fi
\put(-25,5){{\Large ${\cal R}^{\mathbb{R}}$:}}
 \put(0,10){\Large non-}
 \put(0,0){\Large Hermitian}
 \put(45,5){\Large Hermitian}
 \put(105,5){\Large 0}
\put(-25,170){{\Large $\widehat{{\cal R}}{}^{\mathbb{R}}$:}}
 \put(45,170){\Large Hermitian}
 \put(105,170){\Large 0}
 \put(135,155){\large $(n_H+1)$-}
 \put(135,145){\large dimensional}
 \put(135,135){\large conformal}
 \put(135,125){\large hypercomplex}
 \put(145,40){\large $n_H$ -}
 \put(135,30){\large dimensional}
 \put(135,20){\large quaternionic}
\put(64,127){\Large $\widehat{\xi }$ }
 \put(40,45){\Large $\xi $}
\end{picture}
\end{center}
\end{figure}

The figure \ref{fig:mapoverview} gives a schematic overview of our
results. The two blocks represent the families of large (upper block) and
small spaces (lower block), where the horizontal lines indicate how they
are related by $\widehat{\xi }$, resp. $\xi $, transformations. They
connect parametrizations of the same manifold with different complex
structures for hypercomplex manifolds, and different affine and $\SU(2)$
connections for the quaternionic manifolds. On the far right, the spaces
have no $\mathbb{R}$ curvature, and part of these are hyper-K{\"a}hler, resp.
quaternionic-K{\"a}hler. The latter two classes are indicated by the thick
lines. The vertical arrows represent the map described in this review,
connecting the manifolds with similar parametrizations. The map between
hypercomplex and quaternionic manifolds is not one point to one point on
this picture, as a hypercomplex manifold is represented by a full
horizontal line where each point is a particular parametrization. The
same holds for the quaternionic manifolds. The vertical lines are a
representation of the map between the horizontal lines. For some
manifolds there is a representation as a hyper-K{\"a}hler (or
quaternionic-K{\"a}hler) manifold. The thick arrow indicates the map between
hyper-K{\"a}hler and quaternionic-K{\"a}hler spaces.

\section{Symmetries}
 \label{ss:symmetries}

In this section we introduce symmetries as generalizations of isometries.
Indeed, we do not have necessarily a metric and thus no Killing vectors.
We will comment on the moment maps of quaternionic manifolds and show how
they originate from the map between conformal hypercomplex and
quaternionic spaces.

If there is a metric, symmetries are generated from Killing vectors
$k_I^X$, where $I$ labels the different generators. From the Killing
equation one derives that
\begin{equation}
  \nabla_X\nabla_Y \Blue{k_I^Z}=\Purple{R{}_{XWY}{}^Z}\Blue{k_I^W}.
 \label{KillingR}
\end{equation}
This condition is independent of a metric. It turns out that a shift of
the coordinates in the Euler-Lagrange equations (\ref{ELphi}) by an
amount $k_I$ leaves the set of these equations invariant if
(\ref{KillingR}) is satisfied. Hence this becomes the defining equation
for symmetries when there is no good metric available. In the presence of
a quaternionic structure, we also demand that the Lie derivative of the
quaternionic structure rotates them, i.e.
\begin{equation}
  {\cal L}_{\Blue{k_I}} \vec J=\vec r_I\times \vec J,
 \label{LiekIJ}
\end{equation}
for some triplet of functions $\vec{r}_I$. These are called quaternionic
symmetries. For hypercomplex (and hyper-K{\"a}hler) manifolds, the
hypercomplex structures should be invariant, which means that
$\vec{r}_I=0$ in the above equation, and the symmetries are then called
triholomorphic.

It can be shown that this equation is equivalent to the existence of a
decomposition of the derivatives of $k_I$ in an $\SU(2)$ part and a
$\Gl(r,\mathbb{H})$ part, similar to the decomposition of the curvature
in (\ref{curvatrel}):
\begin{equation}
  \nabla_X\Blue{k^Y_I}=\nu \vec J_X{}^Y\cdot \Maroon{\vec P_I} +
L_X{}^Y{}_A{}^B \Red{t_{IB}{}^A}.
 \label{decompdk}
\end{equation}
This defines the triplet `moment map' $\vec P_I$. Its value is related to
the triplet $\vec{r}_I$ in (\ref{LiekIJ}):
\begin{equation}
  \nu \Maroon{\vec{P}_I}\equiv -\ft12 \vec{r}_I-
\OliveGreen{\vec{\omega }}(\Blue{k_I}).
 \label{Pinr}
\end{equation}

Since connections can change by $\xi $-transformations, so will the
symmetry condition depend on the choice of $\xi $. We have shown
\cite{Bergshoeff:2004nf} that a symmetry is preserved only under $\xi$
transformations that satisfy
\begin{equation}
  {\cal L}_{k_I}\xi =0.
 \label{LiekIxi}
\end{equation}
Under such a transformation the moment map transforms as
\begin{equation}
  \nu \Maroon{\widetilde {\vec{P}}_I} = \nu\Maroon{ \vec{P}_I}
  -\Red{\xi}(\vec{J}\Blue{k_I}).
 \label{transfMomentMapxi}
\end{equation}

How do symmetries relate under the map between conformal hypercomplex and
quaternionic manifolds? First it is important to realize that the closed
homothetic Killing vector is itself a triholomorphic symmetry according
to the definition (\ref{KillingR}) and  (\ref{LiekIJ}). For any other
`symmetry' of the conformal hypercomplex manifold we demand that it
commutes with the closed homothetic Killing vector. For such symmetries
there is then a one-to-one mapping with symmetries of the quaternionic
submanifold. The precise relation is expressed by giving the components
of this symmetry vector in the large space:
\begin{equation}
  \Blue{\widehat{k}^0_I} =0,\qquad \Blue{\widehat{k}^\alpha _I}=\Red{\vec{k}^\alpha} \cdot \vec{r}_I,
  \qquad
  \Blue{\widehat{k}^X_I}=\Blue{k_I^X(q)}.
 \label{symmcomponents}
\end{equation}
The last relations says that the symmetry vector of the small space is
just the symmetry vector of the large space projected to the small space.
The components $\widehat{k}^\alpha _I$ that disappear after this
projection are related to the moment maps, in the sense that they define
the $\vec{r}_I$.

\section{Conclusions and final remarks}
 \label{ss:conclusions}

We have illustrated how the picture of quaternionic-like manifolds can be
extended by a mapping as in table \ref{tbl:MapQuatlikeMan}.
\begin{table}[t]
  \caption{\it The 1-to-1 map in quaternionic-like manifolds}
\label{tbl:MapQuatlikeMan}
\begin{center}
  \begin{tabular}{||c|c||}\hline\hline
   \textit{\textbf{CONFORMAL}} & \textit{\textbf{CONFORMAL}}  \\
   \textit{\textbf{hypercomplex}} & \textit{\textbf{hyper-K{\"a}hler}}
   \\[3mm]
   $\updownarrow$ & $\updownarrow$\\[3mm]
    \textit{\textbf{quaternionic}} & \textit{\textbf{quaternionic-K{\"a}hler}} \\
\hline\hline
  \end{tabular}
\end{center}
\end{table}
To obtain this map we started with a conformal hypercomplex manifold of
dimension $4n_H+4$, i.e. a manifold admitting a closed homothetic Killing
vector $k$. This vector defines slices of the manifold $z^0=$constant,
where $z^0$ is the coordinate such that $k$ has only the non-zero
component $k^0$. This slice defines a $4n_H+3$ dimensional manifold which
is a tri-Sasakian manifold. These manifolds still possess an $\SU(2)$
symmetry, which can be divided out such that we are left with a $4n_H$
dimensional manifold, which turns out to be quaternionic. Inversely, any
quaternionic manifold can locally be embedded in such a conformal
hypercomplex manifold of dimension 4 higher.

For quaternionic-K{\"a}hler manifolds, the curvature satisfies a relation
(\ref{quatKahlerEinstein}) depending on a number $\nu $ that sets the
scale. We find here that in this picture, $\nu $ depends on the slice,
i.e. $\nu =-(z^0)^{-1}$. In supergravity the value of $z^0$ is the square
of the Planck mass $M_{\rm Planck}^2$. On the other hand, the sign of
$z^0$ determines the signature of the extra 4 dimensions in uplifting a
quaternionic-K{\"a}hler manifold to a conformal hyper-K{\"a}hler manifold. E.g.
if the quaternionic-K{\"a}hler manifold has a completely positive signature,
and $z^0$ is positive ($\nu $ negative), the signature of the
hyper-K{\"a}hler manifold is
\begin{equation}
  (----++++\cdots ++++).
 \label{signhyperK}
\end{equation}
The construction that we presented can be applied to any signature of the
quaternionic space and $\nu $ positive or negative. This gives then an
arbitrary signature $(4p,4q)$ for the hyper-K{\"a}hler manifold. With $\nu $
negative as in supergravity the quaternionic-K{\"a}hler manifold has negative
scalar curvature and this includes non-compact symmetric spaces. But the
construction can thus easily be applied to compact symmetric spaces as
well.

The construction  of the map that we presented uses heavily the
$\xi$-transformations of connections in the quaternionic manifolds. We
have shown that there are analogous $\hat{\xi }$ transformations in
conformal hypercomplex manifolds. The latter transformations act on the
hypercomplex structure, and are in this respect different from the $\xi $
transformations. But the map that we have constructed is compatible with
both these transformations.

Finally, we have discussed the map between triholomorphic symmetries of
the conformal hypercomplex manifold commuting with dilatations and
quaternionic symmetries of the quaternionic manifolds. It was shown that,
apart form the dilatation symmetry, they relate one-to-one. The moment
maps of quaternionic symmetries are related to components of the symmetry
vector of the hypercomplex manifolds orthogonal to the quaternionic
manifold.

%%%%%%%%%%%%%%%%%%%%%%%%%%%%%%%%
\medskip
\section*{Acknowledgments.}

\noindent We are grateful to S. Cucu, T. de Wit and J. Gheerardyn who
contributed to the original paper related to this review.
% zzz for
%interesting and very useful discussions.
We are grateful to D. Alekseevsky, V. Cort{\'e}s and C. Devchand for many
interesting discussions on quaternionic geometry.

This work is supported in part by the European Community's Human
Potential Programme under contract MRTN-CT-2004-005104 `Constituents,
fundamental forces and symmetries of the universe'. The work of A.V.P. is
supported in part by the FWO - Vlaanderen, project G.0235.05 and by the
Federal Office for Scientific, Technical and Cultural Affairs through the
"Interuniversity Attraction Poles Programme -- Belgian Science Policy"
P5/27.

\newpage

\providecommand{\href}[2]{#2}\begingroup\raggedright\endgroup

%%%%%%%%%%%%%%%%%%%%%%%%%%%%%%%%%%%%%%%%%%%%%%%%%%%%%%%%
%\bibliography{refd5conf}

\begin{thebibliography}{10}

\bibitem{Alekseevsky1975}
D.~V. Alekseevsky, \emph{Classification of quaternionic spaces with a
  transitive solvable group of motions}, Math.\ USSR Izvestija {\bf 9} (1975)
297--339
% .

\bibitem{Bagger:1983tt}
J.~Bagger and E.~Witten, \emph{Matter couplings in $N=2$ supergravity},
Nucl.
  Phys. {\bf B222} (1983)
1
%%CITATION = NUPHA,B222,1;%%.

\bibitem{Alekseevsky:2001if}
D.~V. Alekseevsky, V.~Cort{\'e}s, C.~Devchand  and A.~Van~Proeyen,
\emph{Flows
  on quaternionic-K{\"a}hler and very special real manifolds}, Commun. Math.
  Phys. {\bf 238} (2003) 525--543,
\href{http://www.arXiv.org/abs/hep-th/0109094}{{\tt hep-th/0109094}}
%%CITATION = HEP-TH 0109094;%%.

\bibitem{AM1996}
D.~V. Alekseevsky and S.~Marchiafava, \emph{Quaternionic structures on a
  manifold and subordinated structures}, Ann. Matem. pura appl. (IV) {\bf 171}
  (1996)
205--273
% .

\bibitem{Bergshoeff:2004nf}
E.~Bergshoeff, S.~Cucu, T.~de~Wit, J.~Gheerardyn, S.~Vandoren  and
  A.~Van~Proeyen, \emph{The map between conformal hypercomplex /
  hyper-K{\"a}hler and quaternionic(-K{\"a}hler) geometry},
  \href{http://www.arXiv.org/abs/hep-th/0411209}{{\tt hep-th/0411209}},
to be published in Commun. Math. Phys.
%%CITATION = HEP-TH 0411209;%%.

\bibitem{Swann}
A.~Swann, \emph{HyperK{\"a}hler and quaternionic K{\"a}hler geometry}, Math.
  Ann. {\bf 289} (1991)
421--450
% .

\bibitem{PedersenPS1998}
H.~Pedersen, Y.~S. Poon  and A.~F. Swann, \emph{Hypercomplex structures
  associated to quaternionic manifolds}, Diff. Geom. Appl. {\bf 9} (1998)
273--292
% .

\bibitem{deWit:1985px}
B.~de~Wit, P.~G. Lauwers  and A.~Van~Proeyen, \emph{Lagrangians of $N=2$
  supergravity - matter systems}, Nucl. Phys. {\bf B255} (1985)
569
%%CITATION = NUPHA,B255,569;%%.

\bibitem{Galicki:1992tm}
K.~Galicki, \emph{Geometry of the scalar couplings in $N=2$ supergravity
  models}, Class. Quant. Grav. {\bf 9} (1992)
27--40
%%CITATION = CQGRD,9,27;%%.

\bibitem{deWit:1999fp}
B.~de~Wit, B.~Kleijn  and S.~Vandoren, \emph{Superconformal
hypermultiplets},
  Nucl. Phys. {\bf B568} (2000) 475--502,
\href{http://www.arXiv.org/abs/hep-th/9909228}{{\tt hep-th/9909228}}
%%CITATION = HEP-TH 9909228;%%.

\bibitem{deWit:2001dj}
B.~de~Wit, M.~Ro\v{c}ek  and S.~Vandoren, \emph{Hypermultiplets,
  hyperk{\"a}hler cones and quaternion-K{\"a}hler geometry}, JHEP {\bf 02}
  (2001) 039,
\href{http://www.arXiv.org/abs/hep-th/0101161}{{\tt hep-th/0101161}}
%%CITATION = HEP-TH 0101161;%%.

\bibitem{Bergshoeff:2004kh}
E.~Bergshoeff, S.~Cucu, T.~de~Wit, J.~Gheerardyn, S.~Vandoren  and
  A.~Van~Proeyen, \emph{$N = 2$ supergravity in five dimensions revisited},
  Class. Quant. Grav. {\bf 21} (2004) 3015--3041,
\href{http://www.arXiv.org/abs/hep-th/0403045}{{\tt hep-th/0403045}}
%%CITATION = HEP-TH 0403045;%%.

\bibitem{VanProeyen:2003zj}
A.~Van~Proeyen, \emph{Structure of supergravity theories},
  \href{http://www.arXiv.org/abs/hep-th/0301005}{{\tt hep-th/0301005}},
Publicaciones de la Real Sociedad Matem{\'a}tica Espa{\~n}ola (Publications of
  the Royal Spanish Mathematical Society), Vol. 6, eds. J. Fern{\'a}ndez
  N{\'u}{\~n}ez, W. Garc{\'\i}a Fuertes and A. Vi{\~n}a Escalar, pp.3-32.
%%CITATION = HEP-TH 0301005;%%.

\bibitem{Gunaydin:1984bi}
M.~G{\"u}naydin, G.~Sierra  and P.~K. Townsend, \emph{The geometry of $N=2$
  Maxwell--Einstein supergravity and Jordan algebras}, Nucl. Phys. {\bf B242}
  (1984)
244
%%CITATION = NUPHA,B242,244;%%.

\bibitem{deWit:1992cr}
B.~de~Wit and A.~Van~Proeyen, \emph{Broken sigma model isometries in very
  special geometry}, Phys. Lett. {\bf B293} (1992) 94--99,
\href{http://arXiv.org/abs/hep-th/9207091}{{\tt hep-th/9207091}}
%%CITATION = HEP-TH 9207091;%%.

\bibitem{deWit:1984pk}
B.~de~Wit and A.~Van~Proeyen, \emph{Potentials and symmetries of general
gauged
  $N=2$ supergravity -- Yang-Mills models}, Nucl. Phys. {\bf B245} (1984)
89
%%CITATION = NUPHA,B245,89;%%.

\bibitem{Strominger:1990pd}
A.~Strominger, \emph{Special geometry}, Commun. Math. Phys. {\bf 133}
(1990) 163--180
%%CITATION = CMPHA,133,163;%%.

\bibitem{Castellani:1990zd}
L.~Castellani, R.~D'Auria  and S.~Ferrara, \emph{Special geometry without
  special coordinates}, Class. Quant. Grav. {\bf 7} (1990)
1767--1790
%%CITATION = CQGRD,7,1767;%%.

\bibitem{D'Auria:1991fj}
R.~D'Auria, S.~Ferrara  and P.~Fr{\`e}, \emph{Special and quaternionic
  isometries: General couplings in $N=2$ supergravity and the scalar
  potential}, Nucl. Phys. {\bf B359} (1991)
705--740
%%CITATION = NUPHA,B359,705;%%.

\bibitem{Sierra:1983cc}
G.~Sierra and P.~K. Townsend, \emph{An introduction to $N=2$ rigid
  supersymmetry},
in {\em Supersymmetry and Supergravity 1983}, ed. B. Milewski (World
  Scientific, Singapore, 1983)
% .

\bibitem{Gates:1984py}
J.~Gates, S.~James, \emph{Superspace formulation of new nonlinear sigma
  models}, Nucl. Phys. {\bf B238} (1984)
349
%%CITATION = NUPHA,B238,349;%%.

\bibitem{Alekseevsky:1999ts}
D.~V. Alekseevsky, V.~Cort{\'e}s  and C.~Devchand, \emph{Special complex
  manifolds}, J. Geom. Phys. {\bf 42} (2002) 85--105,
\href{http://www.arXiv.org/abs/math.dg/9910091}{{\tt math.dg/9910091}}
%%CITATION = MATH.DG 9910091;%%.

\bibitem{Craps:1997gp}
B.~Craps, F.~Roose, W.~Troost  and A.~Van~Proeyen, \emph{What is special
  K{\"a}hler geometry?}, Nucl. Phys. {\bf B503} (1997) 565--613,
\href{http://www.arXiv.org/abs/hep-th/9703082}{{\tt hep-th/9703082}}
%%CITATION = HEP-TH 9703082;%%.

\bibitem{VanProeyen:1999ya}
A.~Van~Proeyen, \emph{Special K{\"a}hler geometry},
  \href{http://www.arXiv.org/abs/math.dg/0002122}{{\tt math.dg/0002122}},
in the Proceedings of the meeting on quaternionic structures in
mathematics and
  physics, Roma, September 1999; available on
  \verb+http://www.univie.ac.at/EMIS/proceedings/QSMP99/+
%%CITATION = MATH.DG 0002122;%%.

\bibitem{Freed:1997dp}
D.~S. Freed, \emph{Special K{\"a}hler manifolds}, Commun. Math. Phys. {\bf
203}
  (1999) 31--52,
\href{http://www.arXiv.org/abs/hep-th/9712042}{{\tt hep-th/9712042}}
%%CITATION = CMPHA,203,31;%%.

\bibitem{Gates:1984nk}
J.~Gates, S.~J., C.~M. Hull  and M.~Ro\v{c}ek, \emph{Twisted multiplets
and new
  supersymmetric nonlinear sigma models}, Nucl. Phys. {\bf B248} (1984)
157
%%CITATION = NUPHA,B248,157;%%.

\bibitem{Howe:1987qv}
P.~S. Howe and G.~Papadopoulos, \emph{Ultraviolet behavior of
two-dimensional
  supersymmetric nonlinear sigma models}, Nucl. Phys. {\bf B289} (1987)
264--276
%%CITATION = NUPHA,B289,264;%%.

\bibitem{Bergshoeff:2002qk}
E.~Bergshoeff, S.~Cucu, T.~de~Wit, J.~Gheerardyn, R.~Halbersma,
S.~Vandoren
  and A.~Van~Proeyen, \emph{Superconformal $N = 2$, $D = 5$ matter with and
  without actions}, JHEP {\bf 10} (2002) 045,
\href{http://www.arXiv.org/abs/hep-th/0205230}{{\tt hep-th/0205230}}
%%CITATION = HEP-TH 0205230;%%.

\end{thebibliography}
%%Included for WinEdt Gather Purpose (do not remove the comment line below:
%             %input "C:\localtexmf\bibtex\bib\refd5conf.bib"
%             %input "C:\Program Files\MiKTeX\texmf\bibtex\bib\refd5conf.bib"
%\bibliographystyle{toine}
\end{document}